\documentclass[12pt]{article}	

\newtheorem{theorem}{Theorem}
\newtheorem{theorem-s}{Theorem}[section]
\newtheorem{lemma}{Lemma}[section]
\newtheorem{corollary}{Corollary}[section]
\newtheorem{remark}{Remark}[section]
\newtheorem{example}{Example}[section]

\usepackage{amssymb}\usepackage{amsmath}\usepackage{amsfonts}\usepackage{amscd}

\newcommand{\Ker}{\mathrm{Ker}}

\newcommand{\Pol}{\mathrm{Pol}}
\newcommand{\Sol}{\mathrm{Sol}}
\newcommand{\Fin}{\mathrm{Fin}}
\newcommand{\Int}{\mathrm{Int}}
\newcommand{\sgn}{\mathrm{sgn}}
\newcommand{\SQ}{\mathrm{SQ}}

\newcommand{\N}{\mathbb{N}}
\newcommand{\R}{\mathbb{R}}
\newcommand{\Z}{\mathbb{Z}}
\newcommand{\A}{\mathrm{A}}
\newcommand{\B}{\mathrm{B}}
\newcommand{\E}{\mathrm{E}}
\newcommand{\U}{\mathrm{U}}
\newcommand{\Lim}{\mathrm{L}}
\newcommand{\D}{\Delta}
\renewcommand{\O}{\mathrm{O}}
\renewcommand{\o}{\mathrm{o}}
\renewcommand{\S}{\mathrm{S}}

\title{Qualitative approximation of solutions to difference equations} 

\author{ Janusz Migda\\
Faculty of Mathematics and Computer Science,\\ A. Mickiewicz University,
Umultowska 87, 61-614 Pozna\a'n, Poland\\
email: migda@amu.edu.pl}

\begin{document}\maketitle

\begin{abstract} 
We present a new approach to the theory of asymptotic properties of solutions 
of difference equations. Usually, two sequences $x,y$ are called asymptotically 
equivalent if the sequence $x-y$ is convergent to zero i.e., $x-y\in c_0$, where 
$c_0$ denotes the space of all convergent to zero sequences.
We replace the space $c_0$ by various subspaces of $c_0$. Our approach is 
based on using the iterated remainder operator. Moreover, we use the regional 
topology on the space of all real sequences and the `regional' version of the Schauder 
fixed point theorem. 
\end{abstract}

\section{Introduction} 

Let $\N$, $\R$ denote the set of positive integers and the set of real numbers,  respectively. In this paper we assume that 
\[
m\in\N, \qquad f:\R\to\R, \qquad \sigma:\N\to\N, \qquad \lim\sigma(n)=\infty, 
\]
and consider difference equations of the form 

\begin{equation}\label{E}\tag{E}
\Delta^m x_n=a_nf(x_{\sigma(n)})+b_n
\end{equation}
where \ $a_n,b_n\in\R$.
\smallskip\\
Let $p\in\N$. We say that a sequence $x:\N\to\R$ is a $p$-solution of 
equation \eqref{E} if equality \eqref{E} is satisfied for any $n\geq p$. 
We say that $x$ is a solution if it is a $p$-solution for certain $p\in\N$. 
If $x$ is a $p$-solution for any $p\in\N$, then we say that $x$ is a full solution. 
\smallskip\\ 
In this paper, we present a new approach to the theory of asymptotic properties 
of solutions. The main concept, in our theory, is an asymptotic difference pair. 
The idea of the paper is based on the following observation. If $x$ is a solution of \eqref{E}, $f$ is bounded and the sequence $a$ is `sufficiently  small', then $\D^mx$ 
is close to $b$, and $x$ is close to the set 
\[
\D^{-m}b=\{y\in\R^{\N}: \D^my=b\}.
\]
This means that 
\begin{equation}\label{D-mb+Z}
x\in\D^{-m}b+Z
\end{equation}
where $Z$ is a certain space of `small' sequences. Usually $Z=c_0$ is the space of 
all convergent to zero sequences. In this paper we replace $c_0$ by various subspaces 
of $\R^{\N}$. 
\smallskip\\ 
More precisely, assume that $A$ and $Z$ are linear subspaces of $\R^{\N}$ such that  $A\subset\D^mZ$ and $u\alpha\in A$ for any bounded sequence $u$ and any $\alpha\in A$. 
If $a\in A$ and $x$ is a solution of \eqref{E} such that the sequence 
$u=f\circ x\circ\sigma$ is bounded, then 
\[
\D^mx=au+b\in A+b\subset\D^mZ+b.
\]
Hence $\D^mx=\D^mz+b$ for certain $z\in Z$ and we get $\D^m(x-z)=b$. Therefore 
$x-z\in\D^{-m}b$ and we obtain \eqref{D-mb+Z}. 
\smallskip\\ 
We say that $(A,Z)$ is an asymptotic difference pair of order $m$ (precise definition 
is given in Section 3). In classic case, 
for example in  \cite{Popenda 1990}, \cite{Migda 2001}, \cite{Migda 2010}, we have 
\[
A=\{a\in\R^{\N}: \sum_{n=1}^\infty n^{m-1}|a_n|<\infty\}, \qquad Z=c_0.
\]
In this paper we present some other examples of asymptotic difference pairs. Our 
purpose is to present some basic properties of such pairs. Next, we use asymptotic  difference pairs to the study of asymptotic properties of solutions. For a given 
asymptotic difference pair $(A,Z)$, assuming $a\in A$, we obtain sufficient conditions  under which for any solution $x$ of \eqref{E} there exists $y\in\D^{-m}b$ such that 
$x-y\in Z$. Moreover, assuming $Z\subset c_0$ and using fixed point theorem, we obtain  sufficient conditions  under which for any $y\in\D^{-m}b$ there exists a solution $x$ of  \eqref{E} such that  $y-x\in Z$. Even more, we can `compute modulo $Z$' some parts of 
the set of solutions of \eqref{E} (see Theorem \ref{Th4} and Theorem \ref{Th5} in 
Section 4). 
\smallskip\\
The concept of asymptotic difference pair is an effect of comparing the results from 
some previous papers. In those papers, implicitly, some concrete asymptotic difference  pairs are used (for details see Section 7). In fact, this paper is a continuation of 
a cycle of papers \cite{Migda 2001}-\cite{J Migda 2014 d}. 
\smallskip\\
The paper is organized as follows. In Section 2, we introduce notation and terminology. 
In Section 3, we define asymptotic difference pairs and establish some of their basic  properties. In Section 4, we obtain our main results. 
In Section 5, we present some examples of difference pairs. In our investigations 
the spaces $\A(t)$ (see \eqref{A(t)}) play an important role. In Section 6, we obtain 
some characterizations of $\A(t)$. These results extend some classic tests for 
absolute convergence of series and extend results from \cite{J Migda 2014 c}. 
In Section 7, we present some consequences of our main results. Next we give some 
remarks.

\section{Notation and terminology} 

Let $\Z$ denote the set of all integers. If $p,k\in\Z$, $p\leq k$, then $\N_p$, 
$\N_p^k$ denote the sets defined by 
\[
\N_p=\{p,p+1,\dots\}, \qquad  \N_p^k=\{p,p+1,\dots,k\}.
\]
The space of all sequences $x:\N\to\R$ we denote by $\SQ$. We use the symbols 
\[
\Sol(E), \quad \Sol_p(E), \quad \Sol_{\infty}(E)
\]
to denote the set of all full solutions of \eqref{E}, the set of all $p$-solutions of  \eqref{E}, and the set of all solutions of \eqref{E} respectively. Note that 
\[
\Sol(E)\subset\Sol_p(E)\subset\Sol_{\infty}(E)
\]
for any $p\in\N$. For $p\in\N$ we define 
\[
\Fin(p)=\{x\in\SQ:\ x_n=0 \quad\text{for}\quad n\geq p\}.
\]
Moreover, let  
\[
\Fin(\infty)=\Fin=\bigcup_{p=1}^\infty\Fin(p).
\]
Note that all $\Fin(p)$ are linear subspaces of $\SQ$ and 
\[
0=\Fin(1)\subset\Fin(2)\subset\Fin(3)\subset\dots\subset\Fin(\infty).
\]
If $x,y$ in $\SQ$,  then $xy$ denotes the sequence defined by pointwise multiplication 
\[
xy(n)=x_ny_n. 
\]
Moreover, $|x|$ denotes the sequence defined by $|x|(n)=|x_n|$ for every $n$.

\begin{remark}\label{sol}
A sequence $x\in\SQ$ is a $p$-solution of \eqref{E} if and only if 
\[
\D^mx\in a(f\circ x\circ\sigma)+b+\Fin(p)
\]
and, consequently, $x$ is a solution of \eqref{E} if and only if 
\[
\D^mx\in a(f\circ x\circ\sigma)+b+\Fin.
\]
\end{remark}
We use the symbols `big $\O$' and `small $\o$' in the usual sense but for $a\in\SQ$ we 
also regard $\o(a)$ and $\O(a)$ as subspaces of $\SQ$. More precisely, let 
\[
\o(1)=\{x\in\SQ: \ x \ \text{is convergent to zero}\}, \quad 
\O(1)=\{x\in\SQ: \ x \ \text{is bounded}\}
\]
and for $a\in\SQ$ let 
\[
\o(a)=a\o(1)+\Fin=\{ax: \ x\in\o(1)\}+\Fin, 
\]
\[
\O(a)=a\O(1)+\Fin=\{ax: \ x\in\O(1)\}+\Fin.
\]
Moreover, let
\[
\o(n^{-\infty})=\bigcap_{s\in\R}\o(n^s)=\bigcap_{k=1}^\infty\o(n^{-k}), \qquad 
\O(n^{\infty})=\bigcup_{s\in\R}\O(n^s)=\bigcup_{k=1}^\infty\O(n^k).
\]
Note that if $a_n\neq 0$ for any $n$, then 
\[
\o(a)=a\o(1), \qquad \O(a)=a\O(1).
\]
For $b\in\SQ$ and $X\subset\SQ$ we define 
\[
\D^{-m}b=\{y\in\SQ:\ \D^my=b\}, \qquad \D^{-m}X=\{y\in\SQ:\ \D^my\in X\}.
\]
Moreover, let 
\[
\Pol(m-1)=\D^{-m}0=\Ker\D^m=\{x\in\SQ:\D^mx=0\}.
\]
Then $\Pol(m-1)$ is the space of all polynomial sequences of degree less 
than $m$. 
\smallskip
For a subset $A$ of a metric space $X$ and $\varepsilon>0$ we define 
an $\varepsilon$-framed interior of $A$ 
by 
\[
\Int(A,\varepsilon)=\{x\in X:\overline{\B}(x,\varepsilon)\subset A\}
\]
where $\overline{\B}(x,\varepsilon)$ denotes a closed ball of radius $\varepsilon$ 
about $x$. 
\smallskip\\
We say that a subset $U$ of $X$ is a uniform neighbourhood of a subset 
$Z$ of $X$, if there exists a positive number $\varepsilon$ such that 
$Z\subset\Int(U,\varepsilon)$. For a positive constant $M$ let 
\[
|f\leq M|=\{t\in\R: |f(t)|\leq M\}.
\]
Let 
\[
\A(1):=\{a\in\SQ: \sum_{n=1}^\infty|a_n|<\infty\}.
\]
For $t\in[1,\infty)$ we define 
\begin{equation}\label{A(t)}
\A(t):=\{a\in\SQ: \sum_{n=1}^\infty n^{t-1}|a_n|<\infty\}=(n^{1-t})\A(1).
\end{equation}
Moreover, let 
\[
A(\infty)=\bigcap_{t\in[1,\infty)}\A(t)=\bigcap_{k=1}^\infty\A(k).
\]
Obviously any $\A(t)$ is a linear subspace of $\o(1)$ such that 
\[
\O(1)\A(t)\subset\A(t). 
\]
Note that if $1\leq t\leq s$ then 
\[
\A(\infty)\subset\A(s)\subset\A(t)\subset\A(1).
\]

\begin{remark}
If $p\in\N$, $\lambda\in(0,1)$, $t\in[1,\infty)$, $s\in(0,\infty)$, and $\mu>1$, 
then 
\[
0=\Fin(1)\subset\Fin(p)\subset\Fin\subset\o(\lambda^n)\subset\O(\lambda^n)
\subset\o(n^{-\infty})=\A(\infty),
\]
\[
\A(\infty)\subset\A(t)\subset\A(1)\subset\o(1)\subset\O(1)\subset\o(n^s)
\subset\O(n^s)\subset\O(n^\infty)\subset\o(\mu^n)\subset\O(\mu^n).
\]
\end{remark}

\subsection{Unbounded functions} 

We say that a function $g:\R\to\R$ is unbounded at a point $p\in[-\infty,\infty]$ if 
there exists a sequence $x\in\SQ$ such that $\lim\limits_{n\to\infty}x_n=p$ and 
the sequence $g\circ x$ is unbounded. Let 
\[
U(g)=\{p\in[-\infty,\infty]: g \text{ is unbounded at } p\}.
\]
A function $g:\R\to\R$ is called locally bounded if for any $t\in\R$ there exists a neighbourhood $U$ of $t$ such that the restriction $g|U$ is bounded. Note that any 
continuous function and any monotonic function $g:\R\to\R$ is locally bounded.

\begin{remark}
Assume $g:\R\to\R$. Then  
\begin{enumerate}
\item[$(a)$] $g$ is bounded if and only if $U(g)=\emptyset$, 
\item[$(b)$] $g$ is locally bounded if and only if $U(g)\subset\{\infty,-\infty\}$.
\end{enumerate}
\end{remark}

\begin{example} 
Assume $g:\R\to\R$, $T=\{t_1,t_2,\dots,t_n\}\subset\R$. Then 
\begin{enumerate}
\item[$(01)$] $U(\max(1,t))=U(t+|t|)=U(e^t)=\{\infty\}$, 
\item[$(02)$] $U(\min(1,t))=U(t-|t|)=U(e^{-t})=\{-\infty\}$,
\item[$(03)$] if $g$ is a nonconstant polynomial, then $U(g)=\{-\infty,\infty\}$. 
\item[$(04)$] if $g(t)=1/t$ for $t\neq 0$, then $U(g)=\{0\}$, 
\item[$(05)$] if $g(t)=((t-t_1)\cdots(t-t_n))^{-1}$ for $t\notin T$, then $U(g)=T$.
\end{enumerate}
\end{example}

\begin{remark}
Assume $g,h:\R\to\R$. Then 
\[
U(g+h)\subset U(g)\cup U(h), \qquad U(gh)\subset U(g)\cup U(h).
\]
It follows from the fact that if $g$ and $h$ are bounded at a  point $p$, then $g+h$ 
and $gh$ are also bounded at $p$. Note also that if $U(g)\cap U(h)=\emptyset$, then 
\[
U(g+h)=U(g)\cup U(h).
\]
It is a consequence of the fact that if exactly one of the functions $g,h$ is bounded 
at a point $p$, then $g+h$ is unbounded at $p$. 
\end{remark}

\subsection{Regional topology}

For a sequence $x\in\SQ$ we define a generalized norm  $\|x\|\in[0,\infty]$ by 
\[
\|x\|=\sup\{|x_n|:\ n\in\N\}.
\]
We say that a subset $Q$ of $\SQ$ is ordinary if $\|x-y\|<\infty$ 
for any $x,y\in Q$. We regard every ordinary subset $Q$ of $\SQ$ 
as a metric space with metric defined by 
\begin{equation}\label{dxy}
d(x,y)=\|x-y\|.
\end{equation}
Let $U\subset\SQ$. We say that $U$ is regionally open if $U\cap Q$ 
is open in $Q$ for any ordinary subset $Q$ of $\SQ$. The family of 
all regionally open subsets is a topology on $\SQ$ which we call the regional 
topology. We regard any subset of $\SQ$ as a topological space with topology induced 
by the regional topology. The basic properties of regional topology are presented 
in \cite{J Migda 2014 d}. For the purpose of this paper it is enough to use metrices  defined by \eqref{dxy}. We will use the following `regional' version of the Schauder 
fixed point theorem. 

\begin{lemma}\label{FPL}
Assume $y\in\SQ$, $\rho\in\o(1)$ and 
\[
S=\{x\in\SQ:\ |x-y|\leq|\rho|\}.
\]
Then every continuous map $H:S\to S$ has a fixed point.
\end{lemma}
\textbf{Proof.} See \cite[Lemma 4.7]{J Migda 2014 a}. $\quad\Box$

\subsection{Remainder operator}

Let $\S(m)$ denote the set of all sequences $a\in\SQ$ such that the series 
\[
\sum_{i_1=1}^\infty\sum_{i_2=i_1}^\infty\dots\sum_{i_m=i_{m-1}}^\infty a_{i_m}.
\]
is convergent. For any $a\in\S(m)$ we define the sequence $r^m(a)$ by 
\begin{equation}\label{rman}
r^m(a)(n)=\sum_{i_1=n}^\infty\sum_{i_2=i_1}^\infty\dots\sum_{i_m=i_{m-1}}^\infty a_{i_m}. 
\end{equation}
Then $\S(m)$ is a linear subspace of $\o(1)$, $r^m(a)\in\o(1)$ for any $a\in\S(m)$ 
and 
\[
r^m: \S(m)\to\o(1)
\]
is a linear operator which we call the iterated remainder operator of order $m$. 
The value $r^m(a)(n)$ we denote also by $r^m_n(a)$ or simply $r^m_na$. If 
$|a|\in\S(m)$, then $a\in\S(m)$ and $r^m(a)$ is given by 
\[
r^m(a)(n)=\sum_{j=n}^\infty\binom{m-1+j-n}{m-1}a_j.
\]
Note that if $m=1$, then 
\[
r(a)(n)=r^1(a)(n)=\sum_{j=n}^\infty a_j
\]
is $n$-th remainder of the series $\sum_{j=1}^\infty a_j$. 
The following lemma is a consequence of \cite[Lemma 3.1]{J Migda 2014 c}.

\begin{lemma}\label{reo}
Assume $x,u\in\SQ$ and $p\in\N$. Then 
\begin{enumerate}
\item[$(01)$] $|x|\in\S(m)$ if and only if $x\in\A(m)$, 
\item[$(02)$] $x\in\A(m)$ if and only if $\O(x)\subset\S(m)$, 
\item[$(03)$] if $x\in\A(m)$, $u\in\O(1)$, then $|r^m(ux)|\leq|u|r^m|x|$,
\item[$(04)$] if $ux\in\A(m)$, $\D u\geq 0$ and $u>0$, then $ur^m|x|\leq r^m|ux|$,
\item[$(05)$] if $x\in\A(m)$, then $r^m_p|x|\leq\sum_{n=p}^\infty n^{m-1}|x_n|$, 
\item[$(06)$] $\D^m\o(1)=\S(m)$, \quad $r^m\S(m)=\o(1)$,
\item[$(07)$] if $x\in \S(m)$, then $\D^mr^mx=(-1)^mx$,
\item[$(08)$] if $x\in\o(1)$, then $r^m\D^mx=(-1)^mx$, 
\item[$(09)$] if $Z$ is a linear subspace of $\o(1)$, then $r^m\D^mZ=Z$, 
\item[$(10)$] if $A$ is a linear subspace of $\S(m)$, then $\D^mr^mA=A$,
\item[$(11)$] if $x,y\in \S(m)$ and $x_n\leq y_n$ for $n\geq p$, then 
              $r^m_nx\leq r^m_ny$ for $n\geq p$,
\item[$(12)$] $r^m\Fin(p)=\Fin(p)=\D^m\Fin(p)$, $r^m\Fin=\Fin=\D^m\Fin$,
\item[$(13)$] if $x\geq 0$, then $r^mx$ is nonnegative and nonincreasing.
\end{enumerate}
\end{lemma}

\section{Asymptotic difference pairs}

Let $Z$ be a linear subspace of $\SQ$. We say that a subset $W$ of $\SQ$ is 
$Z$-invariant if $W+Z\subset W$. We say, that a subset $X$ of $\SQ$ is: 

\begin{itemize}
\item[] asymptotic if $X$ is $\Fin$-invariant, 
\item[] evanescent if $X\subset\o(1)$, 
\item[] modular if $\O(1)X\subset X$,
\item[] c-stable if $X$ is $\o(1)$-invariant,
\end{itemize}
We say that a pair $(A,Z)$ of linear subspaces of $\SQ$ is a diffference asymptotic 
pair of order $m$ or, simply, $m$-pair if $Z$ is asymptotic, $A$ is modular and 
$A\subset\D^mZ$. We say that an $m$-pair $(A,Z)$ is evanescent if $Z$ is evanescent. 

\begin{remark}\label{dp-rem1} 
For any $a\in\SQ$ the spaces $\o(a)$ and $\O(a)$ are asymptotic and modular.
\end{remark}

\begin{remark}\label{dp-rem2} 
If $W$ is an asymptotic subset of $\SQ$, $x\in W$ and $x'$ is a sequence obtained 
from $x$ by change finite number of terms, then $x'\in W$. Moreover, a linear subspace 
$Z$ of $\SQ$ is asymptotic if and only if $\Fin\subset Z$. 
\end{remark}

\begin{remark}\label{dp-rem3} 
If $(A,Z)$ is an evanescent $m$-pair, then, using Lemma \ref{reo} (06), we have 
\[
A\subset\D^mZ\subset\D^m\o(1)=\S(m)\subset\o(1).
\]
Hence the space $A$ is evanescent.
\end{remark}

\begin{remark}\label{dp-rem4} 
If $a\in\SQ$, then the sequence $\sgn\circ a$ is bounded and $|a|=(\sgn\circ a)a$. 
Hence, if  $W$ is a modular subset of $\SQ$, then $|a|\in W$ for any $a\in W$. In  particular, if $(A,Z)$ is an evanescent $m$-pair and $a\in A$, then 
\[
|a|\in A\subset\D^mZ\subset\D^m\o(1)=\S(m).
\]
Therefore $A\subset\A(m)$ and, for any $a\in A$, the sequences $r^ma$ and $r^m|a|$ 
are defined.
\end{remark}

\begin{lemma}\label{L1}
Assume $(A,Z)$ is an $m$-pair, $a,b\in\SQ$, and $a-b\in A$. Then 
\[
\D^{-m}a+Z=\D^{-m}b+Z.
\]
\end{lemma}
\textbf{Proof.} 
We have $a-b\in A\subset\D^mZ$. Hence there exists $z_0\in Z$ such that  
$a-b=\D^mz_0$. Let $x\in\D^{-m}a$ and $z\in Z$. Then  
\[
\D^m(x-z_0)=\D^mx-\D^mz_0=a-(a-b)=b. 
\] 
Therefore $x+z=x-z_0+z_0+z\in\D^{-m}b+Z$. Thus 
\[
\D^{-m}a+Z\subset\D^{-m}b+Z.
\]
Since $b-a=-(a-b)\in A$, we have 
\[
\D^{-m}b+Z\subset\D^{-m}a+Z. 
\]
The proof is complete. $\quad\Box$

\begin{lemma}\label{L2}
Assume $(A,Z)$ is an $m$-pair and $b\in A$. Then 
\[
\D^{-m}b+Z=\Pol(m-1)+Z. 
\]
\end{lemma}
\textbf{Proof.} 
This lemma is an immediate consequence of the previous lemma. $\quad\Box$ 

\begin{lemma}\label{L3}
Assume $(A,Z)$ is an $m$-pair, $a\in A$, $b,x\in\SQ$ and 
\[ 
\D^mx\in\O(a)+b
\]
Then $x\in\D^{-m}b+Z$.
\end{lemma}
\textbf{Proof.} 
The condition $a\in A$ implies $\O(a)\subset A$. Hence, 
\[
\D^mx-b\in\O(a)\subset A\subset\D^mZ.
\] 
Therefore, there exists $z\in Z$ such that $\D^mx-b=\D^mz$. Then 
\[
\D^m(x-z)=\D^mx-\D^mz=b.
\]
Thus $x-z\in\D^{-m}b$ and we obtain $x=x-z+z\in\D^{-m}b+Z$. $\quad\Box$

\begin{lemma}[Comparison test] 
Assume $A$ is an asymptotic, modular linear subspace of $\SQ$, $b\in A$, $a\in\SQ$, 
and $|a_n|\leq|b_n|$ for large $n$. Then $a\in A$. 
\end{lemma}
\textbf{Proof.} 
Assume $|a_n|\leq|b_n|$ for $n\geq p$. Let
\[
h_n=
\begin{cases}
0 & \text{if\; $b_n=0$ \; }\\
a_n/b_n & \text{if \quad $b_n\neq 0$. \quad }
\end{cases}
\]
Then $h\in\O(1)$. Moreover, if $n\geq p$ and $b_n=0$, then $a_n=0$. Hence 
$a_n=h_nb_n$ for $n\geq p$. Therefore $a-hb\in\Fin(p)$. Let $z=a-hb$. Then 
\[
a=hb+z\in\O(1)A+\Fin\subset A+A=A.
\]

\section{Solutions} 
In this section, in Theorems \ref{Th4} and \ref{Th5} we obtain our main results. 
First we introduce the notion of $f$-ordinary and $f$-regular sets. We use these sets in  Theorem \ref{Th5}. At the end of the section we present some examples of  $f$-regular sets. 
\smallskip\\
We say that a subset $W$ of $\SQ$ is $f$-ordinary if for any $x\in W$ the sequence 
$f\circ x$ is bounded. We say that a subset $W$ of $\SQ$ is $f$-regular if for any 
$x\in W$ there exists an index $p$ such that $f$ is continuous and bounded on some 
uniform neighborhood of the set $x(\N_p)$. For $x\in\SQ$ let 
\[
\Lim(x)=\{p\in[-\infty,\infty]: p \text{ is a limit point of } x\}. 
\]

\begin{lemma}\label{LB1} 
If $x\in\SQ$, then 
\[
f\circ x\in\O(1) \quad\text{or}\quad \Lim(x)\cap\U(f)\neq\emptyset.
\] 
\end{lemma}
\textbf{Proof.} 
Assume the sequence $f\circ x$ is ubounded from above. Then there exists a subsequence 
$x_{n_k}$ such that 
\[
\lim_{k\to\infty}f(x_{n_k})=\infty.
\]
Let $y_k=x_{n_k}$ and let $p\in L(y)$. There exists a subsequence $y_{k_i}$ such that 
\[
\lim_{i\to\infty}y_{k_i}=p.
\]
Then $\lim_{i\to\infty}f(y_{k_i})=\infty$ and we obtain $p\in U(f)$. Since $y$ is a 
subsequence of $x$, we have $\Lim(y)\subset\Lim(x)$. Hence $p\in U(f)\cap L(x)$. 
Analogously, if the sequence $f\circ x$ is unbounded from below, then  
$U(f)\cap\Lim(x)\neq\emptyset$. $\quad\Box$ 
\smallskip\\
Note that if the sequence $f\circ x$ is bounded, then $f\circ x\circ\sigma$ is 
also bounded.

\begin{theorem}\label{LU}
Assume $(A,Z)$ is an $m$-pair, $a\in A$, and $x\in\Sol_{\infty}(\E)$. Then 
\[
x\in\D^{-m}b+Z \quad\text{or}\quad \Lim(x)\cap\U(f)\neq\emptyset.
\]
\end{theorem}
\textbf{Proof.} 
Assume $\Lim(x)\cap\U(f)=\emptyset$. Then, by Lemma \ref{LB1}, the sequence 
$f\circ x$ is bounded. Hence the sequence $f\circ x\circ\sigma$ is bounded too. 
By Remark \ref{sol}, 
\[
\D^mx\in a(f\circ x\circ\sigma)+b+\Fin.
\]
Hence 
\[
\D^mx\in a\O(1)+\Fin+b=\O(a)+b.
\]
Using Lemma \ref{L3} we obtain $x\in\D^{-m}b+Z$. The proof is complete. $\quad\Box$ 

\begin{corollary} 
Assume $(A,Z)$ is an $m$-pair, $a\in A$, $B\cup C=\R$, $C$ is closed in $\R$, $f$ 
is bounded on $B$, $\U(f)\subset\R$, and $x$ is a solution of \eqref{E}. Then 
\[
x\in\D^{-m}b+Z \quad\text{or}\quad \Lim(x)\cap C\neq\emptyset.
\] 
\end{corollary}
\textbf{Proof.} 
Using the relation $\U(f)\subset\R$ we see that $\U(f)\subset C$. Hence  
the assertion is a consequence of Theorem \ref{LU}. $\quad\Box$

\begin{example}
Assume $(A,Z)$ is an $m$-pair, $a\in A$, and $f(t)=t^{-1}$ for $t\neq 0$. If $x$ is 
a solution of \eqref{E} such that $0$ is not a limit point of $x$, then, by 
Theorem \ref{LU}, $x\in\D^{-m}b+Z$.
\end{example}

\begin{example}
Assume $(A,Z)$ is an $m$-pair, $a\in A$, $f$ is continuous and there exists a 
proper limit $\lim_{t\to\infty}f(t)$. Then, by Theorem \ref{LU}, for any bounded 
below solution $x$ of \eqref{E} we have $x\in\D^{-m}b+Z$. 
\end{example}

\begin{theorem}\label{abc}
Assume $(A,Z)$ is an $m$-pair, $a\in A$, and $W\subset\SQ$ is $f$-ordinary. Then 
\[
W\cap\Sol_{\infty}(\E)\subset\D^{-m}b+Z.
\]
\end{theorem}
\textbf{Proof.} 
Let $x\in W\cap\Sol_{\infty}(\E)$. By Remark \ref{sol}, 
\[
\D^mx\in a(f\circ x\circ\sigma)+b+\Fin.
\]
Since $x\in W$, we have $f\circ x=\O(1)$. Hence $f\circ x\circ\sigma=\O(1)$ and 
\[
\D^mx\in a\O(1)+\Fin+b=\O(a)+b.
\]
Now, the assertion follows from Lemma \ref{L3}. $\quad\Box$

\begin{theorem}\label{Th3} 
Assume $(A,Z)$ is an evanescent $m$-pair, $a\in A$, $M>0$, $R=Mr^m|a|$, $f$ is continuous  on $|f\leq M|$, $y\in\D^{-m}b$, $p\in\N$, and 
\[
(y\circ\sigma)(\N_p)\subset\Int(|f\leq M|,R_p). 
\]
Then $y\in\Sol_p(\E)+Z$.
\end{theorem} 
\textbf{Proof.} 
For $x\in\SQ$ let 
\[
x^*=f\circ x\circ\sigma.
\]
Let $\rho_n=R_n$ for $n\geq p$ and $\rho_n=0$ for $n<p$. Moreover, let 
\[
T=\{x\in\mathrm{SQ}:|x-y|\leq R_p\}, \qquad S=\{x\in\mathrm{SQ}:|x-y|\leq\rho\}.
\]
Obviously $S\subset T$. Let $x\in S$. If $n\geq p$, then 
$|x_{\sigma(n)}-y_{\sigma(n)}|\leq R_p$. Hence 
\[
x_{\sigma(n)}\in\overline{\B}(y_{\sigma(n)},R_p)\subset|f\leq M|.
\]
Therefore for $x\in S$ we have $|x^*|\leq M$ and 
\begin{equation}\label{ax*}
ax^*\in A\subset\mathrm{S}(m). 
\end{equation}
Let 
\begin{equation}\label{H}
H:S\to\mathrm{SQ}, \qquad 
H(x)(n)=\left\{
\begin{array}{lll}
y_n & \text{for} & n<p\\
y_n+(-1)^mr^m_n(ax^*) & \text{for} & n\geq p.
\end{array}
\right.
\end{equation}
If $x\in S$ and $n\geq p$, then 
\[
|H(x)(n)-y_n|=|r^m_n(ax^*)|\leq r^m_n|ax^*|\leq Mr^m_n|a|=R_n=\rho_n. 
\] 
Hence $HS\subset S$. As in the proof of \cite[Theorem 1]{Migda 2010} one can show that 
$H$ is continuous and there exists a sequence $x\in S$ such that $Hx=x$. Then 
\begin{equation}\label{xnyn}
x_n=y_n+(-1)^mr^m_n(ax^*)
\end{equation}
for $n\geq p$. Therefore, for $n\geq p$ we have 
\[
\D^mx_n=b_n+a_nf(x_{\sigma(n)}).
\]
Thus $x\in\Sol_p(\E)$. By \eqref{xnyn} 
\[
y-x+(-1)^mr^m_n(ax^*)\in\Fin(p)
\]
Hence, by \eqref{ax*}
\begin{equation}\label{y-x}
y-x\in r^mA+\Fin.
\end{equation}
Using definition of evanescent $m$-pair and Lemma \ref{reo} (06) we have 
\[
A\subset\D^mZ\subset\D^m\o(1)=\S(m).
\]
Hence, using Lemma \ref{reo} (09), we have 
\begin{equation}\label{rmA}
r^mA\subset r^m\D^mZ=Z.
\end{equation}
Now, by \eqref{y-x}, we obtain 
\[
y-x\in Z+\Fin=Z.
\]
Hence $y\in x+Z$. The proof is complete. $\quad\Box$

\begin{corollary}\label{cor1} 
Assume $(A,Z)$ is an evanescent $m$-pair, $a\in A$, $y\in\D^{-m}b$ and $\{y\}$ is 
$f$-regular. Then 
\[
y\in\Sol_{\infty}(\E)+Z.
\]
\end{corollary}
\textbf{Proof.} 
There exist a positive $M$ and $\delta>0$ such that 
\[
(y\circ\sigma)(\N)\subset\Int(|f\leq M|,\delta). 
\]
Let $R=Mr^m|a|$. Then $R=\o(1)$ and $R_p<\delta$ for certain $p$. Hence  
\[
\Int(|f\leq M|,\delta)\subset\Int(|f\leq M|,R_p)
\]
and, by Theorem \ref{Th3}, $y\in\Sol_p(\E)+Z$. $\quad\Box$ \\
The next theorem is our first main result. We assume that $f$ is continuous and 
bounded. This assumption is very strong but our result is also strong. 

\begin{theorem}\label{Th4} 
Assume $(A,Z)$ is an evanescent $m$-pair, $a\in A$, $p\in\N$, 
and $f$ is continuous and bounded. Then 
\[
\Sol(\E)+Z=\Sol_p(\E)+Z=\Sol_{\infty}(\E)+Z=\D^{-m}b+Z.
\]
\end{theorem}
\textbf{Proof.} 
Choose $M$ such that $|f|\leq M$. Then $|f\leq M|=\R$. Hence 
\[
\Int(|f\leq M|,\delta)=\R 
\]
for any positive $\delta$. By Theorem \ref{Th3} we have 
\[
\D^{-m}b\subset\Sol_p(\E)+Z
\]
for any $p$. For a given $p\in\N$ we obtain 
\[
\D^{-m}b+Z\subset\Sol(\E)+Z\subset\Sol_p(\E)+Z\subset
\Sol_{\infty}(\E)+Z.
\]
On the other hand, by Theorem \ref{abc}, taking $W=\SQ$ we obtain 
\[
\Sol_{\infty}(\E)+Z\subset\D^{-m}b+Z.
\]
The proof is complete. $\quad\Box$

\begin{lemma}\label{AL} 
Assume $Z$ is a linear subspace of a linear space $X$, $D,S,W\subset X$, 
$W$ is $Z$-invariant, $W\cap S\subset D+Z$ and $W\cap D\subset S+Z$. Then 
\[
W\cap S+Z=W\cap D+Z.
\]
\end{lemma}
\textbf{Proof.} 
Assume $w\in W\cap S$. Since $W\cap S\subset D+Z$, we have $w\in W\cap(D+Z)$.  
Hence, there exist $d\in D$ and $z\in Z$ such that $w=d+z$. Since $W$ is 
$Z$-invariant, we obtain 
\[
d=w-z\in W+Z\subset W.
\]
Hence $w=d+z\in(W\cap D)+Z$. Therefore $W\cap S\subset W\cap D+Z$ and we obtain 
\[
W\cap S+Z\subset W\cap D+Z+Z=W\cap D+Z.
\]
Analogously, we obtain $W\cap D+Z\subset W\cap S+Z$. $\quad\Box$ \\
Now we are ready to prove our second main result.

\begin{theorem}\label{Th5} 
Assume $(A,Z)$ is an evanescent $m$-pair, $a\in A$, and $W\subset\SQ$. Then 
\begin{itemize}
\item[$(a)$] if $W$ is $f$-ordinary, then $W\cap\Sol_{\infty}(\E)\subset\D^{-m}b+Z$, 
\item[$(b)$] if $W$ is $f$-regular, then $W\cap\D^{-m}b\subset\Sol_{\infty}(\E)+Z$, 
\item[$(c)$] if $W$ is $f$-regular and $Z$-invariant, then 
\[
W\cap\Sol_{\infty}(\E)+Z=W\cap\D^{-m}b+Z.
\]
\end{itemize}
\end{theorem} 
\textbf{Proof.} 
Assertion (a) is a special case of Theorem \ref{abc}. (b) is a consequence of 
Corollary \ref{cor1}. Using (a), (b), Lemma \ref{AL} and the fact that any 
$f$-regular set $W\subset\SQ$ is also $f$-ordinary we obtain (c). $\quad\Box$ 

\begin{remark}\label{inv-rem1}
Any subset of an $f$-regular set is $f$-regular. If $Z$ is a linear subspace of 
$\o(1)$, then any c-stable subset $W$ of $\SQ$ is also $Z$-invariant. 
\end{remark} 

\begin{remark}\label{inv-rem2} 
Assume $W\subset\SQ$ is $f$-regular and $Z$ is a linear subspace of $\o(1)$. 
Then the set $W+Z$ is $f$-regular and $Z$-invariant. 
\end{remark}

\begin{example}\label{ex1} 
Assume $f$ is continuous and bounded on a certain uniform neighborhood of a set 
$Y\subset\R$. Then the set 
\[
W=\{y\in\SQ:\, y(\N)\subset Y\}
\]
is $f$-regular. If $x\in\SQ$ and $z\in\o(1)$, then $\Lim(x+z)=\Lim(x)$. Hence the sets 
\[
W_1=\{y\in\SQ:\, \Lim(y)\subset Y\}, \qquad W_2=\{y\in\SQ:\, \lim y_n\in Y\} 
\]
are $f$-regular and c-stable. 
\end{example}

\begin{example}\label{ex2} 
If $f$ is bounded, then $\SQ$ is $f$-ordinary and  c-stable. Moreover, if $f$ 
is continuous, then $\SQ$ is $f$-regular. 
\end{example} 

\begin{example}\label{ex3} 
If $f$ is locally bounded, then the set $\O(1)$ of all bounded sequences is 
$f$-ordinary and c-stable. Moreover, if $f$ is continuous, then $\O(1)$ 
is $f$-regular. 
\end{example} 

\begin{example}\label{ex4} 
If $f$ is locally bounded, then the set $C$ of all convergent sequences 
is $f$-ordinary and c-stable. Moreover, if $f$ is continuous, then $C$ is $f$-regular. 
\end{example} 

\begin{example}\label{ex5} 
Let $Z$ be a linear subspace of $\o(1)$ and $p\in\N$. We say that a sequence $x\in\SQ$ 
is $(p,Z)$-asymptotically periodic if there exists a $p$-periodic sequence $y$ 
such that $x-y\in Z$. If $f$ is locally bounded, then the set $W$ of all 
$(p,Z)$-asymptotically periodic sequences is $f$-ordinary and $Z$-invariant. 
Moreover, if $f$ is continuous, then $W$ is $f$-regular. 
\end{example} 

\begin{example}
If $f(t)=e^t$, then the sets 
\[
W_1=\{x\in\SQ: \limsup_{n\to\infty}x_n<\infty\} \quad\text{and}\quad 
W_2=\{x\in\SQ: \lim_{n\to\infty}x_n=-\infty\}
\]
are $f$-regular and c-stable. 
\end{example}

\begin{example}
If $f$ is continuous and $\limsup\limits_{t\to\infty}|f(t)|<\infty$, then 
the sets  
\[
W_1=\{x\in\SQ: \liminf_{n\to\infty}x_n>-\infty\} \quad\text{and}\quad 
W_2=\{x\in\SQ: \lim_{n\to\infty}x_n=\infty\}
\]
are $f$-regular and c-stable. 
\end{example}

\begin{example}
If $f(t)=t^{-1}$ for $t\neq 0$, then the set $W=\{x\in\SQ: 0\notin\Lim(x)\}$ 
is $f$-regular and c-stable. 
\end{example}

\begin{example}
Assume $g:\R\to\R$ is continuous, $T=\{t_1,t_2,\dots,t_n\}\subset\R$ and 
\[
f(t)=\frac{g(t)}{(t-t_1)(t-t_2)\dots(t-t_n)}
\]
for $t\notin T$. Then the set $W=\{x\in\SQ: T\cap\Lim(x)=\emptyset\}$ 
is $f$-regular and c-stable. 
\end{example}

\section{Examples of difference pairs} 

We say that a subset $A$ of $\SQ$ is an $m$-space, if $(A,A)$ is an $m$-pair. 
In this section we present some examples of difference $m$-pairs and $m$-spaces. 
Next we establish some lemmas to justify our examples.

\begin{remark} 
Assume that $(A,Z)$ is an $m$-pair. If $Z^*$ is a linear subspace of $\SQ$ such 
that $Z\subset Z^*$, then $(A,Z^*)$ is an $m$-pair. Analogously, if $A_*$ is a 
modular subspace of $A$, then $(A_*,Z)$ is an $m$-pair. 
\end{remark}

\begin{example}\label{mp1}
If $a\in\A(m)$, then $(\O(a), r^m\O(a))$ is an evanescent $m$-pair.
\end{example}

\begin{example}\label{mp2}
If $X$ is an asymptotic and modular subspace of $\A(m)$, then $(X, r^mX)$ is an 
evanescent $m$-pair. 
\end{example}

\begin{example}\label{mp3}
Let $s\in(-\infty,-m)$. The following pairs are evanescent $m$-pairs 
\[ 
(\o(n^{s}),\, \o(n^{s+m})), \qquad (\O(n^{s}),\, \O(n^{s+m})).
\]
\end{example} 

\begin{example}\label{mp4}
If $s\in\R$ and $(s+1)(s+2)\dots(s+m)\neq 0$, then 
\[ 
(\o(n^{s}),\, \o(n^{s+m})), \qquad (\O(n^{s}),\, \O(n^{s+m})) 
\]
are $m$-pairs.
\end{example} 

\begin{example}\label{mp5}
Let $\lambda\in(0,1)$. The following spaces are evanescent $m$-spaces 
\[
\Fin, \quad \o(\lambda^n), \quad \O(\lambda^n), \quad \o(n^{-\infty}). 
\]
\end{example} 

\begin{example}\label{mp6}
Let $\lambda\in(1,\infty)$. The following spaces are $m$-spaces 
\[
\o(\lambda^n), \quad \O(\lambda^n), \quad \O(n^{\infty}). 
\]
\end{example} 

\begin{example}\label{mp7}
If $s\in(-\infty,0]$, then $(\A(m-s),\,\o(n^s))$ is an evanescent $m$-pair. 
\end{example}

\begin{example}\label{mp8}
Assume $s\in(-\infty,m-1]$, and $q\in\N_0^{m-1}$. Then 
\[
(\A(m-s),\,\o(n^s)), \qquad (\A(m-q),\, \D^{-q}\o(1))
\]
are $m$-pairs.
\end{example}

\begin{example}\label{mp9}
If $t\in[1,\infty)$, then $(\A(m+t),\ \A(t))$ is an evanescent $m$-pair.
\end{example}
Note that Example \ref{mp1} is a special case of Example \ref{mp2}. Note also that 
Lemma \ref{reo} (10) justifies Example \ref{mp2}. To justify Examples \ref{mp3} and 
\ref{mp4} we need the following four lemmas. 

\begin{lemma}[Cesaro-Stolz lemma]\label{DLH}
Assume $x,y\in\SQ$, $y$ is strictly monotonic and one of the following conditions 
is satisfied
\begin{itemize}
\item[$(a)$] $x=\o(1)$ and $y=\o(1)$, 
\item[$(b)$] $y$ is unbounded. 
\end{itemize}
Then
\[
\liminf\frac{\D x}{\D y}\leq \liminf\frac{x}{y}\leq
\limsup\frac{x}{y}\leq \limsup\frac{\D x}{\D y}.
\]
\end{lemma}
\textbf{Proof.} 
If (a) is satisfied, then the assertion is proved in \cite{Agarwal}. Assume (b) and 
$y$ is unbounded from above. Then $y$ is increasing and $\lim y_n=\infty$. Let
\[
L=\liminf\frac{\Delta x_n}{\Delta y_n}.
\]
If $L=-\infty$, then the inequality 
\begin{equation}\label{dlh-e1}
\liminf\frac{\D x}{\D y}\leq \liminf\frac{x}{y}
\end{equation}
is obvious. Assume $L>-\infty$. Choose a constant $M$ such that $M<L$. Then there 
exists an index $p$ such that $\Delta x_n/\Delta y_n\geq M$ for $n\geq p$. We can 
assume that $y_n>0$ and $\Delta y_n>0$ for $n\geq p$. If $n\geq p$, then 
\[
x_n-x_p=\Delta x_p+\Delta x_{p+1}+\dots+\Delta x_{n-1}
\]
\[
\geq M(\Delta y_p+\Delta y_{p+1}+\dots+\Delta y_{n-1})=M(y_n-y_p).
\]
Hence $x_n\geq My_n+x_p-My_p$ and 
\[
\frac{x_n}{y_n}\geq M+\frac{x_p-My_p}{y_n}
\]
for $n\geq p$. Since $\lim(1/y_n)=0$, we have
\[
\liminf\frac{x_n}{y_n}\geq M.
\]
Therefore, we obtain \eqref{dlh-e1}. Similarly, one can prove the inequality 
\[
\limsup\frac{x}{y}\leq \limsup\frac{\D x}{\D y}.
\]
Replacing $y$ by $-y$ we obtain the result if $y$ is unbounded from below. $\quad\Box$

\begin{lemma}\label{DLH3}
Assume $x\in \SQ$, $s\in\R$, and $s>-1$ or $x=\o(1)$. Then 
\[
\D x=\o(n^s)\ \Longrightarrow\ x=\o(n^{s+1}), \qquad
\D x=\O(n^s)\ \Longrightarrow\ x=\O(n^{s+1}).
\]
\end{lemma}
\textbf{Proof.} 
If $s=-1$, then, by assumption, $x=\o(1)=\o(n^{s+1})$. Hence the assertion is true 
for $s=-1$. Assume $s\neq -1$. Note that 
\[
\frac{\D x_n}{\D n^{s+1}}=\frac{\D x_n}{n^s}\frac{n^s}{\D n^{s+1}}
\]
By the proof of \cite[Lemma 2.1]{Migda 2013}, the sequence $(n^s/\D n^{s+1})$ is  convergent. Hence the assertion follows from Lemma \ref{DLH}. $\quad\Box$

\begin{lemma}\label{DLH4}
Assume $s\in\R$ and $s+1\neq 0$. Then 
\[   
\o(n^s)\subset\D\o(n^{s+1}), \qquad \O(n^s)\subset\D\O(n^{s+1}).
\]
\end{lemma}
\textbf{Proof.} 
Assume $z=\o(n^s)$. Choose $x\in \SQ$ such that $z=\D x$. If $s>-1$, then, by 
Lemma \ref{DLH3}, $x=\o(n^{s+1})$. Let $s<-1$. Then the series $\sum z_n$ is 
convergent. Let
\[
\sigma=\sum_{n=1}^\infty z_n, \quad x_1=0, \quad x_n=z_1+\dots+z_{n-1}-\sigma \quad 
\text{for} \quad  n>1.
\]
Then $x=\o(1)$, $\D x=z$ and by Lemma \ref{DLH3}, we have $x=\o(n^{s+1})$. Hence we
obtain $\o(n^s)\subset\D\o(n^{s+1})$. Analogously $\O(n^s)\subset\D\O(n^{s+1})$. 
$\quad\Box$ 

\begin{lemma}\label{DLH5}
Assume $s\in\R$ and $(s+1)(s+2)\dots(s+m)\neq 0$. Then 
\[
\o(n^s)\ \subset\ \D^m\o(n^{s+m}), \qquad \O(n^s)\ \subset\ \D^m\O(n^{s+m}).
\]
\end{lemma}
\textbf{Proof.} 
The assertion is an easy consequence of the previous lemma. $\quad\Box$ 
\medskip\\
Lemma \ref{DLH5} justify Examples \ref{mp3} and \ref{mp4}.

\begin{lemma}\label{lambda} 
If $\lambda\in(0,1)\cup(1,\infty)$, then 
\[ 
\o(\lambda^n)\subset\D^m\o(\lambda^n), \qquad 
\O(\lambda^n)\subset\D^m\O(\lambda^n).
\]
\end{lemma}
\textbf{Proof.} 
Let $x,w\in\SQ$ and $\D w=x$. Since $\D\lambda^n=\lambda^{n+1}-\lambda^n=\lambda^n(\lambda-1)$, we have 
\begin{equation}\label{geo}
\frac{\D w_n}{\D\lambda^n}=\frac{x_n}{\D\lambda^n}=L\frac{x_n}{\lambda^n}
\end{equation}
where, $L=1/(\lambda-1)$. Assume $\lambda\in(0,1)$ and $x\in\o(\lambda^n)$. Then 
the series $\sum_{n=1}^\infty x_n$ is convergent. 
Hence $x\in\S(1)=\D\o(1)$ and there exists $w\in\o(1)$ such that $x=\D w$. Using 
\eqref{geo} and the fact that $x\in\o(\lambda^n)$ we have $\D w_n/\D\lambda^n\to 0$. 
Moreover, $w_n\to 0$ and $\lambda^n\to 0$. 
By Lemma \ref{DLH}, we obtain $w\in\o(\lambda^n)$. Hence 
\[
x=\D w\in\D\o(\lambda^n). 
\] 
Therefore $\o(\lambda^n)\subset\D\o(\lambda^n)$ and, by induction, 
\[
\o(\lambda^n)\subset\D^m\o(\lambda^n).
\] 
If $x\in\O(\lambda^n)$, then the sequence $x_n/\lambda^n$ is bounded and, by \eqref{geo},  the sequence $\D w_n/\D\lambda^n$ is also bounded. Hence, by Lemma \ref{DLH},  $w\in\O(\lambda^n)$ and we obtain $\O(\lambda^n)\subset\D\O(\lambda^n)$. Moreover, 
by induction 
\[
\O(\lambda^n)\subset\D^m\O(\lambda^n). 
\] 
If $\lambda>1$, then $\lambda^n\to\infty$ and using Lemma \ref{DLH} (b) we obtain 
the result. $\quad\Box$

\begin{lemma}\label{o-inf} 
$\o(n^{-\infty})\subset\D^m\o(n^{-\infty})$.
\end{lemma}
\textbf{Proof.} 
Using \cite[Lemma 4.8]{J Migda 2014 c} we have $r\A(k+1)\subset\A(k)$ for any 
$k\in\N$. Hence $r\A(\infty)\subset r\A(k+1)\subset\A(k)$ and we get 
\[
r\A(\infty)\subset\bigcap_{k\in\N}\A(k)=\A(\infty).
\]
Therefore $r^2\A(\infty)=rr\A(\infty)\subset r\A(\infty)\subset\A(\infty)$ and so on. 
After $m$ steps we obtain 
\[
r^m\A(\infty)\subset\A(\infty).
\]
By Lemma \ref{reo} (07), we have $\D^mr^m\A(\infty)=\A(\infty)$. Hence 
\[
\A(\infty)=\D^mr^m\A(\infty)\subset\D^m\A(\infty).
\]
Now the result follows from the equality $\o(n^{-\infty})=\A(\infty)$. $\quad\Box$ 
\medskip\\
Using Lemma \ref{reo} (12), Lemma \ref{lambda} and Lemma \ref{o-inf} we justify 
Example \ref{mp5}. By Lemma \ref{DLH5}, we have 
$\O(n^s)\subset\D^m\O(n^\infty)$ for any $s>m$. Hence 
\[
\O(n^\infty)=\bigcup_{s>m}\O(n^s)\subset\D^m\O(n^\infty).
\]
Therefore, using Lemma \ref{lambda}, we obtain Example \ref{mp6}.

\begin{lemma}\label{Am-s} 
If $s\in(-\infty,m-1]$, then $\A(m-s)\subset\D^m(\o(n^s))$.
\end{lemma}
\textbf{Proof.} 
Let $a\in\A(m-s)$. Choose $x\in\SQ$ such that $a=\D^mx$. By 
\cite[Theorem 2.1]{Migda 2013} we have $x\in\Pol(m-1)+\o(n^s)$. Hence 
\[
a=\D^mx\in\D^m(\Pol(m-1)+\o(n^s))
\]
\[
=\D^m\Pol(m-1)+\D^m\o(n^s)=\D^m\o(n^s).
\]

\begin{lemma}\label{Am-q}
If $q\in\N_0^{m-1}$, then $\A(m-q)\subset\D^m\D^{-q}\o(1)$.
\end{lemma}
\textbf{Proof.} 
Let $a\in\A(m-q)$. Choose $x\in\SQ$ such that $a=\D^mx$. By 
\cite[Lemma 3.1 (d)]{J Migda 2014 a} we have $x\in\Pol(m-1)+\D^{-q}\o(1)$. Hence 
\[
a=\D^mx\in\D^m(\Pol(m-1)+\D^{-q}\o(1))=\D^m\D^{-q}\o(1).
\]
Using Lemmas \ref{Am-s} and \ref{Am-q} we obtain Examples \ref{mp7} and \ref{mp8}.

\begin{lemma}\label{Ex9}
If $t\in[m+1,\infty)$, then $r^m\A(t)\subset\A(t-m)$.
\end{lemma}
\textbf{Proof.} 
Choose $k\in\N$ such that $k\leq t<k+1$. Let $s=t-k$. Then  
\[
A(t)=n^{1-t}\A(1)=n^{1-(k+s)}\A(1)=n^{-s}n^{1-k}\A(1)=n^{-s}\A(k).
\]
Hence, for $a\in\A(t)$ we have $n^sa\in\A(k)$. 
By Lemma \ref{reo} (04), 
\[
a\in\A(k) \quad\text{and}\quad n^sr|a|\leq r|n^sa|. 
\]
Since $|n^sa|\in\A(k)$ and  $r(\A(k))\subset\A(k-1)$, we have 
\[
r|n^sa|\in\A(k-1). 
\]
By the comparison test we obtain $n^sr|a|\in\A(k-1)$. 
Using the inequality $|ra|\leq r|a|$ we have $n^s|ra|\leq n^sr|a|$. By comparison 
test, $n^s|ra|\in\A(k-1)$. Hence 
\[
ra\in n^{-s}\A(k-1)=\A(t-1). 
\]
Therefore 
\[ 
r(\A(t))\subset\A(t-1)
\]
and, by induction, we obtain the result. $\quad\Box$ 
\medskip\\
Now let $t\in[1,\infty)$. By Lemma \ref{Ex9} we have $r^m\A(m+t)\subset\A(t)$. 
Hence, using Lemma \ref{reo} (07), 
\[
\A(m+t)=\D^mr^m\A(m+t)\subset\D^m\A(t)
\]
and we obtain Example \ref{mp9}.

\section{Absolute summable sequences} 

In our investigations the spaces $\A(t)$ play an important role. In this section 
we obtain some characterizations of $\A(t)$. Our results extend some classical tests 
for absolute convergence of series and extend results from \cite{J Migda 2014 c}. 

\begin{lemma} 
Assume $t\in[1,\infty)$ and $s\in\R$. Then $(n^s)\in\A(t)\Leftrightarrow s<-t$.
\end{lemma}
\textbf{Proof.} 
We have 
\[
(n^s)\in\A(t)\Leftrightarrow (n^s)\in(n^{1-t})\A(1)\Leftrightarrow 
(n^{t+s-1})\in\A(1)\Leftrightarrow t+s-1<-1 \Leftrightarrow s<-t.
\]

\begin{lemma}[Generalized logarithmic test]\label{test1}
Assume $a\in\SQ$, $t\in[1,\infty)$ and 
\[
u_n=-\frac{\ln|a_n|}{\ln n}.
\]
Then
\begin{enumerate}
\item[$(1)$] \ if \ $\liminf u_n>t$, \ then \ $a\in\A(t)$,
\item[$(2)$] \ if \ $u_n\leq t$ \ for large \ $n$, \ then \ $a\notin\A(t)$,
\item[$(3)$] \ if \ $\limsup u_n<t$, \ then \ $a\notin\A(t)$,
\item[$(4)$] \ if \ $\lim u_n=\infty$, \ then \ $a\in\A(\infty)$.
\end{enumerate}
\end{lemma}
\textbf{Proof.} 
If $\liminf u_n>t$, then there exists a number $s>t$ such that $u_n>s$ for large $n$. 
Then $|a_n|\leq n^{-s}$ for large $n$. Hence (1) follows from the comparison test 
and from the fact that $(n^{-s})\in\A(t)$. If $u_n\leq t$ for large $n$, then 
$|a_n|\geq n^{-t}$ for large $n$. Hence (2) follows from the fact that  $(n^{-t})\notin\A(t)$. The assertion (3) follows immediately from (2) and (4) is a  consequence of (1). $\quad\Box$

\begin{lemma}[Generalized Raabe's test]\label{Raabe}\label{test2}
Assume \ $a\in\SQ$, \ $t\in[1,\infty)$, 
\[
u_n=n\left(\frac{|a_n|}{|a_{n+1}|}-1\right). \qquad \text{Then}
\]
\begin{enumerate}
\item[$(1)$] \ if \ $\liminf u_n>t$, \ then \ $a\in\A(t)$,
\item[$(2)$] \ if \ $u_n\leq t$ \ for large \ $n$, \ then \ $a\notin\A(t)$,
\item[$(3)$] \ if \ $\limsup u_n<t$, \ then \ $a\notin\A(t)$,
\item[$(4)$] \ if \ $\lim u_n=\infty$, \ then \ $a\in\A(\infty)$.
\end{enumerate}
\end{lemma}
\textbf{Proof.} 
Let
\[
b_n=n^{t-1}a_n, \qquad w_n=n\left(\frac{|b_n|}{|b_{n+1}|}-1\right). \qquad 
\text{Then}
\]
\[
w_n=n\left(\frac{n^{t-1}|a_n|}{(n+1)^{t-1}|a_{n+1}|}-1\right)=
n\left(\left(\frac{n}{n+1}\right)^{t-1}
\frac{|a_n|}{|a_{n+1}|}-1\right)
\]
\[
=n\left(\frac{n}{n+1}\right)^{t-1}
\left(\frac{|a_n|}{|a_{n+1}|}-\left(\frac{n+1}{n}\right)^{t-1}\right)=
n\left(\frac{n}{n+1}\right)^{t-1}
\left(\frac{|a_n|}{|a_{n+1}|}-\left(1+\frac{1}{n}\right)^{t-1}\right).
\]
If \ $s\in\R$, \ then using the Taylor expansion of the function $(1+x)^s$ we obtain  
\[
(1+x)^s=1+sx+\o(x) \qquad \text{for} \qquad x\to 0. 
\]
\[
\text{Hence} \qquad 
\left(1+\frac{1}{n}\right)^{t-1}=1+(t-1)\frac{1}{n}+\o(n^{-1}). 
\]
Therefore 
\[
w_n=\left(\frac{n}{n+1}\right)^{t-1}n\left(\frac{|a_n|}{|a_{n+1}|}-1\right)-
\left(\frac{n}{n+1}\right)^{t-1}(t-1-n\o(n^{-1}))
\]
\[
=c_nu_n-c_n(t-1-\o(1)), \qquad c_n=\left(\frac{n}{n+1}\right)^{t-1}\to 1.
\]
Thus
\[
\liminf w_n=\liminf u_n-(t-1)=\liminf u_n-t+1.
\]
Hence, if $\liminf u_n>t$, then $\liminf w_n>1$ and by the usual Raabe's test we 
obtain $b\in\A(1)$ i.e., $a\in\A(t)$. The assertion (1) is proved.  
Now, we assume that $u_n\leq t$ for large $n$. Then  
\[
n\left(\frac{|a_n|}{|a_{n+1}|}-1\right)\leq t \quad \text{i.e.,} \quad 
\frac{|a_n|}{|a_{n+1}|}\leq\frac{t}{n}+1 \quad \text{for large} \quad n. 
\]
\[
\text{Hence} \qquad 
w_n=n\left(\frac{n}{n+1}\right)^{t-1}
\left(\frac{|a_n|}{|a_{n+1}|}-\left(1+\frac{1}{n}\right)^{t-1}\right)
\]
\[
\leq n\left(\frac{n}{n+1}\right)^{t-1}
\left(\frac{t}{n}+1-\left(1+\frac{1}{n}\right)^{t-1}\right)
\]
It is easy to see that if $t\geq 1$ and $x\in(0,1)$, then $(1+x)^t\geq 1+tx$. 
\[
\text{Hence} \qquad 
\left(1+\frac{1}{n}\right)^t\geq 1+\frac{t}{n}, \qquad \text{Hence} \qquad 
\frac{t}{n}+1-\left(1+\frac{1}{n}\right)\left(1+\frac{1}{n}\right)^{t-1}\leq 0. 
\]
\[
\text{Hence} \qquad 
\frac{t}{n}+1-\left(1+\frac{1}{n}\right)^{t-1}\leq 
\frac{1}{n}\left(1+\frac{1}{n}\right)^{t-1}=
\frac{1}{n}\left(\frac{n+1}{n}\right)^{t-1}
\]
Hence $w_n\leq 1$ for large $n$ and, by the usual Raabe's test, we obtain  
$b\notin\A(1)$ i.e., $a\notin\A(t)$. The assertion (2) is proved. 
(3) is an immediate consequence of (2). (4) follows from (1). $\quad\Box$

\begin{lemma}[Generalized Schl\"{o}milch's test]\label{test3}
Assume $a\in\SQ$, $t\in[1,\infty)$, 
\[
u_n=n\ln\frac{|a_n|}{|a_{n+1}|}. 
\]
Then
\begin{enumerate}
\item[$(1)$] if $\liminf u_n>t$, then $a\in\A(t)$,
\item[$(2)$] if $u_n\leq t$ for large $n$, then $a\notin\A(t)$,
\item[$(3)$] if $\limsup u_n<t$, then $a\notin\A(t)$,
\item[$(4)$] if $\lim u_n=\infty$, then $a\in\A(\infty)$.
\end{enumerate}
\end{lemma}
\textbf{Proof.} 
If $\liminf u_n=b>t$ and $c\in(t,b)$, then $\liminf u_n>c$ for large $n$. Hence  
\[
\frac{|a_n|}{|a_{n+1}|}\geq\exp\left(\frac{c}{n}\right). 
\] 
Since $e^x\geq 1+x$ for $x>0$, we have  
\[
\frac{|a_n|}{|a_{n+1}|}>1+\frac{c}{n} \qquad \text{and} \qquad 
n\left(\frac{|a_n|}{|a_{n+1}|}-1\right)>c>t
\]
for large $n$. Now, by Raabe's test we obtain (1). \\ 
Assume $u_n\leq t$ for large $n$. Then 
\[
\ln\frac{|a_n|}{|a_{n+1}|}\leq\frac{t}{n} \qquad\text{and}\qquad 
\frac{|a_n|}{|a_{n+1}|}\leq e^{\frac{t}{n}}
\]
for large $n$. Let $b_n=(n-1)^{-t}$. Since 
\[
e<\left(1+\frac{1}{n-1}\right)^n, 
\]
we have 
\[
e^{\frac{t}{n}}<\left(1+\frac{1}{n-1}\right)^t=\left(\frac{n}{n-1}\right)^t=
\frac{b_n}{b_{n+1}}.
\]
Hence 
\[
\frac{|a_n|}{|a_{n+1}|}\leq e^{\frac{t}{n}}<\frac{b_n}{b_{n+1}} 
\qquad\text{and}\qquad \frac{|a_n|}{b_n}<\frac{|a_{n+1}|}{b_{n+1}} 
\]
for large $n$. Hence, there exists a $\lambda>0$ such that $|a_n|/b_n>\lambda$ 
for large $n$. Therefore
\[
|a_n|>\lambda b_n>\lambda n^{-t}
\]
for large $n$. Using the fact that $(n^{-t})\notin\A(t)$ we have $a\notin\A(t)$ 
and we obtain (2). The assertion (3) is an immediate consequence of (2). (4) follows 
from (1). $\quad\Box$

\begin{lemma}[Generalized Gauss's test]\label{test4}
Let $a\in\SQ$, $t\in[1,\infty)$, $\lambda, s\in\R$, $s<-1$ and 
\[
\frac{|a_n|}{|a_{n+1}|}=1+\frac{\lambda}{n}+\O(n^s)
\]
Then 
\begin{enumerate}
\item[$(a)$] if $\lambda>t$, \ then \ $a\in\A(t)$, 
\item[$(b)$] if $\lambda\leq t$, \ then \ $a\notin\A(t)$.
\end{enumerate}
\end{lemma}
\textbf{Proof.} 
Let $b_n=n^{t-1}a_n$. Then  
\[
\frac{|b_n|}{|b_{n+1}|}=\frac{n^{t-1}|a_n|}{(n+1)^{t-1}|a_{n+1}|}=
\left(\frac{n}{n+1}\right)^{t-1}\frac{|a_n|}{|a_{n+1}|}
\]
\[
=\left(\frac{n}{n+1}\right)^{t-1}\left(1+\frac{\lambda}{n}+\O(n^s)\right).
\]
It is easy to see that  
\[
\left(\frac{n}{n+1}\right)^{t-1}=\left(1-\frac{1}{n+1}\right)^{t-1}=
1-\frac{t-1}{n+1}+\O\left(\frac{1}{(n+1)^2}\right)=
1-\frac{t-1}{n}+\O\left(\frac{1}{n^2}\right)
\]
Hence
\[
\frac{|b_n|}{|b_{n+1}|}=\left(1-\frac{t-1}{n}+\O(n^{-2})\right)
\left(1+\frac{\lambda}{n}+\O(n^s)\right)
\]
\[
=1+\frac{\lambda}{n}+\O(n^s)-\frac{t-1}{n}-\frac{\lambda(t-1)}{n^2}+
\O(n^{s-1})+\O(n^{-2})
\]
\[
=1+\frac{\lambda-t+1}{n}+\O(n^{s'}).
\]
For some $s'<-1$. If $\lambda>t$, then $\lambda-t+1>1$ and, by the usual Gauss's test, 
$b\in\A(1)$. Hence $a\in\A(t)$. Analogously, if $\lambda\leq t$, then 
$\lambda-t+1\leq 1$ and $b\notin\A(1)$. Therefore $a\notin\A(t)$. $\quad\Box$

\begin{lemma}[Generalized Kummer's test]\label{test5}
Assume $a$, $c$ are positive sequences,  
\[
t\in[1,\infty), \qquad 
K_n=\frac{c_na_n}{a_{n+1}}\left(\frac{n}{n+1}\right)^{t-1}-c_{n+1}. 
\]
Then
\begin{enumerate}
\item[$(1)$] if $\liminf K_n>0$, then $a\in\A(t)$, 
\item[$(2)$] if the series $\sum_{n=1}^\infty c_n^{-1}$ is divergent and $K_n\leq 0$ 
             for large $n$, then $a\notin\A(t)$, 
\item[$(3)$] if the series $\sum_{n=1}^\infty c_n^{-1}$ is divergent and 
             $\limsup K_n<0$, then $a\notin\A(t)$. 
\end{enumerate}
\end{lemma}
\textbf{Proof.} 
This lemma is an easy consequence of the usual Kummer's test since, by definition of 
the space $\A(t)$, we have 
\[
(a_n)\in\A(t)\ \Leftrightarrow\ (n^{t-1})(a_n)\in\A(1).
\]

\begin{lemma}[Generalized Bertrand's test]\label{test6}
Assume $a\in\SQ$, $t\in[1,\infty)$ and  
\[
\frac{|a_n|}{|a_{n+1}|}=1+\frac{t}{n}+\frac{\lambda_n}{n\ln n}. 
\]
Then 
\begin{enumerate}
\item[$(1)$] if $\liminf\lambda_n>1$, then $a\in\A(t)$, 
\item[$(2)$] if $\lambda_n\leq 1$ for large $n$, then $a\notin\A(t)$, 
\item[$(3)$] if $\limsup\lambda_n<1$, then $a\notin\A(t)$. 
\end{enumerate}
\end{lemma}
\textbf{Proof.} 
Let 
\[
c_n=n\ln n, \quad u_n=\left(\frac{n}{n+1}\right)^{t-1}, \quad 
K_n=u_n\frac{c_na_n}{a_{n+1}}-c_{n+1}. 
\]
Then 
\[
K_n=u_n\left(1+\frac{t}{n}+\frac{\lambda_n}{n\ln n}\right)n\ln n-c_{n+1}=
u_n(n+t)\ln n+u_n\lambda_n-c_{n+1}.
\]
Since 
\[
u_n=\left(\frac{n}{n+1}\right)^{t-1}=\left(1-\frac{1}{n+1}\right)^{t-1}=
1+\frac{1-t}{n+1}+\O(n^{-2}), 
\]
we have 
\[
K_n=\left(1+\frac{1-t}{n+1}+\O(n^{-2})\right)(n+t)\ln n-(n+1)\ln(n+1)+\lambda_n u_n
\]
\[
=(n+t)\ln n+\frac{(n+t)(1-t)\ln n}{n+1}+\o(1)-(n+1)\ln(n+1)+\lambda_n u_n
\]
\[
=(n+1)\ln n+\frac{(n+1)(t-1)\ln n}{n+1}+
\frac{(n+t)(1-t)\ln n}{n+1}+\o(1)-(n+1)\ln(n+1)+\lambda_n u_n
\]
\[
=\ln\left(\frac{n}{n+1}\right)^{n+1}+\lambda_n u_n+\o(1)=
\ln\left(1-\frac{1}{n+1}\right)^{n+1}+\lambda_n u_n+\o(1)
\]
\[
=-1+\lambda_n u_n+\o(1).
\]
Since $u_n\to 1$, we have 
\[
\liminf K_n=-1+\liminf\lambda_n. 
\]
Hence, by Kummer's test and by the divergence of the series 
$\sum_{n=1}^\infty c_n^{-1}$, we obtain (1) and (3). 
Since $(1+x)^t\geq 1+tx$ for $t,x\in[0,\infty)$, we have 
\[
1+\frac{t}{n}\leq\left(1+\frac{1}{n}\right)^t \qquad\text{and}\qquad 
n+t\leq n\left(\frac{n+1}{n}\right)^t. 
\]
Hence 
\[
u_n(n+t)\leq n\frac{n+1}{n}\left(\frac{n+1}{n}\right)^{t-1}u_n=n+1. 
\]
Moreover,
\[
e<\left(1+\frac{1}{n}\right)^{n+1}=\left(\frac{n+1}{n}\right)^{n+1} 
\quad \Rightarrow \quad e^{-1}>\left(\frac{n}{n+1}\right)^{n+1} 
\quad \Rightarrow \quad \ln\left(\frac{n}{n+1}\right)^{n+1}<-1.
\]
Hence 
\[
K_n=u_n(n+t)\ln n+u_n\lambda_n-c_{n+1}\leq
(n+1)\ln n+u_n\lambda_n-(n+1)\ln(n+1)
\]
\[
=u_n\lambda_n+\ln\left(\frac{n}{n+1}\right)^{n+1}<-1+u_n\lambda_n.
\]
Now, using Kummer's test and the fact that $u_n\in(0,1]$, we obtain (2). $\quad\Box$

\section{Remarks} 

In this section we present some consequences of our results. 
Next we give some final remarks. 
\smallskip\\
The first part of Theorem \ref{Th5} we may state in the following form. 

\begin{theorem-s}\label{Th6} 
Assume $(A,Z)$ is an evanescent $m$-pair, $a\in A$, and $W\subset\SQ$. Then 
\begin{itemize}
\item[$(a)$] if $W$ is $f$-ordinary, then for any solution $x$ of \eqref{E} such that 
             $x\in W$ there exists a sequence $y$ such that $\D^my=b$ and $x-y\in Z$, 
\item[$(b)$] if $W$ is $f$-regular, then for any sequence $y\in W$ such that $\D^my=b$ 
             there exists a solution $x$ of \eqref{E} such that $y-x\in Z$. 
\end{itemize}
\end{theorem-s} 
Using this theorem and Lemma \ref{L2} we obtain 

\begin{theorem-s}\label{Th7} 
Assume $(A,Z)$ is an evanescent $m$-pair, $a,b\in A$, and $W\subset\SQ$. Then 
\begin{itemize}
\item[$(a)$] if $W$ is $f$-ordinary, then for any solution $x$ of \eqref{E} such that 
             $x\in W$ there exists a polynomial sequence $\varphi\in\Pol(m-1)$ such 
             that $x-\varphi\in Z$, 
\item[$(b)$] if $W$ is $f$-regular, then for any polynomial sequence 
             $\varphi\in\Pol(m-1)$ such that $\varphi\in W$ there exists a solution 
             $x$ of \eqref{E} such that $\varphi-x\in Z$. 
\end{itemize}
\end{theorem-s} 
Using Theorem \ref{Th6}, Example \ref{mp9} and the generalized Raabe's test 
(Lemma \ref{test2}) we obtain 

\begin{theorem-s}\label{Th8} 
Assume $W\subset\SQ$ is $f$-regular, $t\in[1,\infty)$, and 
\[
\liminf n\left(\frac{|a_n|}{|a_{n+1}|}-1\right)>m+t. 
\]
Then for any $y\in W\cap\D^{-m}b$ there exists a solution $x$ of \eqref{E} such that 
\[
\limsup n\left(\frac{|y_n-x_n|}{|y_{n+1}-x_{n+1}|}-1\right)\geq t. 
\]
\end{theorem-s} 
Using Example \ref{mp7} and the generalized Schl\"{o}milch's test (Lemma \ref{test3}) 
we obtain 

\begin{theorem-s}\label{Th9} 
Assume $W\subset\SQ$ is $f$-ordinary, $s\in(-\infty,0]$, and 
\[
\liminf n\ln\frac{|a_n|}{|a_{n+1}|}>m-s. 
\]
Then for any solution $x$ of \eqref{E} such that $x\in W$ there exists $y\in\D^{-m}b$ 
such that 
\[
x_n-y_n=\o(n^s). 
\]
\end{theorem-s} 
Using Example \ref{mp5} we obtain 

\begin{theorem-s}\label{Th10} 
Assume $W\subset\SQ$ is $f$-regular, $\lambda\in(0,1)$, and $a\in\o(\lambda^n)$.
Then for any $y\in W\cap\D^{-m}b$ there exists a solution $x$ of \eqref{E} such that 
\[
y_n-x_n=\o(\lambda^n). 
\]
\end{theorem-s} 
Using other examples of $m$-pairs one can obtain many other theorems. 
\smallskip\\
Asymptotic difference pairs are used, implicitly, in some papers. 
The classical case is $(\A(m),\o(1))$, see for example \cite{Popenda 1990}, 
\cite{Migda 2001}, \cite{Migda 2010}. 
The pair $(\A(m-s),\o(n^s))$, for a fixed $s\in(-\infty,0]$, is used in 
\cite{Migda 2013}, \cite{J Migda 2014 a}, \cite{J Migda 2014 b} and 
\cite{J Migda 2014 d}. 
The pair $(\A(m+p),\A(p))$, for a fixed $p\in\N$, is used in \cite{J Migda 2014 c}. 
The pair $(\A(m-q),\D^{-q}\o(1))$, for a fixed $q\in\N_0^{m-1}$, is used in 
\cite[Theorem 7.5]{J Migda 2014 a}. 
\smallskip\\
Our results may be partially extended to the case of nonautonomous equations. 
The basic difference is as follows. If $f:\R\to\R$, $x\in\SQ$ and the sequence 
$(f(x_n))$ is bounded, then the sequence $(f(x_{\sigma(n)}))$ is also bounded. 
On the other hand, if $f:\N\times\R\to\R$, $x\in\SQ$, then the boundedness of 
the sequence $(f(n,x_n))$ does not imply the boundedness of the sequence 
$(f(n,x_{\sigma(n)}))$.
\smallskip\\
In some papers the term generalized solution is used instead of our solution and 
the term solution in place of our full solutions.
\smallskip\\
The terminology is a matter of taste. A separate question is why study 
the generalized solution at all. As a kind of motivation we give three examples. 
These examples are taken from \cite{J Migda 2014 b}.

\begin{example}
Assume $a_n\geq 0$, the series $\sum_{n=1}^\infty a_n$ is convergent and there exists 
an index $p>1$ such that $a_p=0$, $a_{p+1}=1$. Consider the equation 
\[
\Delta x_n=a_n|x_{n-1}|. 
\]
Then every number $\lambda\in\R$ is the limit of a certain solution. 
On the other hand, if $x$ is a full convergent solution, then $\lim x_n\geq 0$. 
\end{example}

\begin{example}
Assume $a_n\geq 0$, the series $\sum_{n=1}^\infty na_n$ is convergent and $a_p=1$ for  certain $p$. Consider the equation 
\[
\Delta^2 x_n=a_n x_n^2.
\]
Then every real constant $\lambda$ is the limit of a certain solution but if  $\lambda$ 
is the limit of a full solution, then $\lambda<2$. 
\end{example}

\begin{example}
Assume that the series $\sum_{n=1}^\infty a_n$ is absolutely convergent. 
Consider the equation 
\[
\Delta x_n=a_nx_n+a_n.
\]
Then any real constant $\lambda$ is the limit of a certain solution. Morever, 
if $a_n\neq-1$ for all $n\in\N$, then any real $\lambda$ is the limit of a certain full solution. On the other hand, if $a_p=-1$ for certain $p$ and $x$ is a $p$-solution, 
then $x_n=-1$ for any $n>p$. 
\end{example}

\end{document}